\documentclass{article}

\usepackage{amsmath}
\usepackage{amssymb}
\usepackage{amsthm}
\usepackage{epsfig}

\newtheorem{theorem}{Theorem}

\newtheorem{problem}{Problem}

\newcommand{\bE}{{\bf E}}
\newcommand{\bP}{{\bf P}}
\newcommand{\cG}{{\mathcal G}}
\newcommand{\floor}[1]{\ensuremath{\left \lfloor {#1} \right \rfloor}}
\newcommand{\ceil}[1]{\ensuremath{\left \lceil {#1} \right \rceil}}

\DeclareMathOperator{\sgn}{sgn}

\begin{document}
\title{Where Do Power Laws Come From?}
\author{Joshua Cooper
\thanks{University of South Carolina, Columbia, SC 29208.
{\tt cooper@math.sc.edu}}
\and
Linyuan Lu \thanks{University of South Carolina, Columbia, SC
29208.
{\tt lu@math.sc.edu}}}

\maketitle

\begin{abstract}
    What distribution of graphical degree sequence is invariant under
``scaling''? Are these graphs always power-law graphs? We show the
answer is a surprising ``yes'' for sparse graphs if we ignore  isolated
vertices, or more generally, the vertices with degrees less than a fixed
constant $k$. We obtain a concentration result on the degree sequence
of a random induced subgraph.  The case of hypergraphs (or set-systems) is also  examined.
\end{abstract}
\section{Introduction}


Quite a few recent papers use the term ``scale-free networks'' to
refer to those large sparse graphs formed from real-world data.  Such
graphs often exhibit power-law degree distributions.  Namely, the
number of vertices with degree $d$ is roughly proportional to
$d^{-\beta}$, for some positive $\beta$.  However, the term
``scale-free'' is rarely defined in the literature, at least in the rigorous mathematical sense.  Furthermore, accounts in the literature of how power laws arise have
been largely model-dependent.  That is, a number of models of random-graph
growth have been proposed that give rise, under circumstances of
varying generality, to power-law degree distributions. 
The most popular growth model of this kind is the ``preferential attachment'' scheme, exemplified by \cite{acle,ba,baj,klein}.

 Though many of
the growth rules are quite intuitive -- in that one expects many
real-world phenomena to approximate them -- an explanation of the
sheer ubiquity of power laws that does not appeal to particular models
is conspicuously lacking.


Here we attempt to address these omissions.  First, it is natural to ask, what is a scale?  An obvious candidate for a scale is the number of
vertices of a graph.  Here ``scaling the graph down'' means ``taking
an induced subgraph''.  Of course, subgraphs may look quite different from one another.  Hence, we consider only the average behavior.\\

\noindent
{\bf Random induced subgraph $G_p$:}
  For any $0<p<1$, let $G_p$ be the induced subgraph of $G$ on
a random subset of vertices $S$. For each vertex $v$ of $G$,
$v$ is in $V(G_p)$ with probability $p$ independently.\\

There are some simple cases that the graph $G_p$ is similar to
$G$. For example,
\begin{itemize}
\item Let $G$ be a complete graph on $n$ vertices.
Then $G_p$ is also a complete graph on around $pn$ vertices.
\item Let $G$ be an empty graph on $n$ vertices.
Then $G_p$ is also an empty graph on around $pn$ vertices.
\item For any constant $q\in (0,1)$,
let $G$ be the random graph $G(n, q)$. Then $G_p$ is
also a random graph $G(m,q)$ over a randomly chosen set
of size $m \sim pn$.
\end{itemize}
Crucially, these examples are not ``real-world graphs'', in the sense that graphs appearing ``in nature'' tend to be quite sparse. Most vertices have small degrees. To characterize this property, we use the following definition:\\

\noindent For a given sequence $\{\lambda_d\}_{d=0}^\infty$ satisfying
$\sum_{d=0}^\infty\lambda_d =1$, with $\lambda_d \geq 0$ for all $d \geq 0$, a sequence of graphs $\{G^n\}$ on $n$ vertices
is said to have degree sequence with limit distribution $\{\lambda_d\}_{d=0}^\infty$ if the number of vertices with degree $d$ in $G^n$ is $\lambda_d n +o(n)$ for each $d \geq 0$.  We also say that $\{G^n\}$ has limit distribution $\{\lambda_d\}_{d=k}^\infty$, for $\sum_{d \geq k} \lambda_d \leq 1$, if $G^n$ has $\lambda_d n +o(n)$ vertices of degree $d$ for each $d \geq k$.


We consider two questions.
\begin{enumerate}
\item  If the degree sequence of $G$ in $\{G_n\}$ has a limit
  distribution, for any fixed $p$, does the degree sequence of the
  random induced subgraph $G_p$ also have a limit distribution?
\item For what distribution $\{\lambda_k\}_{k=0}^\infty$ is
the limit distribution of the degree sequence of $G_p$ essentially
the same as the limit distribution of the degree sequence of $G$?
\end{enumerate}

To answer the first question, we observe that a vertex of degree
$cn$ in $G$ would badly affect the concentration of the degree sequence of
$G_p$. On the other hand, using the vertex-exposure martingale,
we can show that the degree sequence of $G_p$ will have a limit
distribution if
$$\sum_v \deg^2(v) = O(n^{2-\epsilon}).$$
This condition is satisfied, for example, if $G$ has maximum degree bounded by $n^{1/2-\epsilon}$.

Suppose $a_0, a_1, a_2, \ldots, $ is the degree frequency sequence of
a graph $G$, with $a_d$ representing the number of vertices in
$G$ with degree $d$. What is the degree frequency sequence of $G_p$?
If a vertex $v$ survives in $G_p$, its degree has binomial distribution $B(d_G(v),p)$. There is no simple way
to describe the joint distribution because of edge-correlations.
Nonetheless, the expected degree frequency sequence for $G_p$ is easy to compute.
Let $b_0, b_1, b_2, \ldots$ be the expected degree frequency sequence of
the random induced subgraph $G_p$. We have
\begin{equation}
  \label{eq:1}
  b_d=p\sum_{k\geq i} a_k{k\choose d} p^d (1-p)^{k-d}
\end{equation}
for all $d=0,1,2,\ldots$.
Note that $\{b_d\}_{d\geq 0}$ depends linearly on
$\{a_d\}_{d\geq 0}$.  We can therefore normalize both sequences by dividing by $n$.

Therefore, from now on, we assume
$a_i$ are the fraction of numbers of vertices with degree $i$
in graph $G$. More precisely, we consider a sequence of graphs
$G_n$, such that the number of vertices with degree $d$ in $G_n$ is
$a_d n+o(n)$.
We only consider sparse graphs such that
\begin{equation}
  \label{eq:2}
  \sum_{i\geq 0} a_i =1.
\end{equation}

It is worth remarking that this manuscript can be read, in effect, as a response to the well-known Stumpf, Wiuf, and May paper, ``Subnets of scale-free networks are not scale-free: Sampling properties of networks'' (\cite{stumpf}) and its authors' related publications.  Although the present authors became aware of this work only after discovering the results below, it is clear that there is a very strong resemblance to the work of Stumpf, Wiuf, and May.  However, we offer the counter-assertion ``Subnets of scale-free networks {\it are} scale-free, as long as one ignores suitably small-degree vertices.''  We also take a somewhat different tack by studying, in particular, the asymptotic conditions under which scale-freeness holds.

\section{Scale-Free Degree Sequences}

Let $A(x)=\sum_{i=0}^\infty a_i x^i$ be the generating function of $\{a_i\}_{i\geq 0}$ and $B(x)=\sum_{i=0}^\infty b_i x^i$ be the generating function of
$\{b_i\}_{i\geq 0}$. Both $A(x)$ and $B(x)$ converge on $[-1,1]$.

We have
\begin{eqnarray*}
  B(x) &=& \sum_{i=0}^n b_i x^i \\
&=& \sum_{i=0}^\infty p\sum_{k\geq i} (a_k + o(1)) {k\choose i} p^i (1-p)^{k-i} x^i\\
&=& p\sum_{k=0}^\infty a_k \sum_{i=0}^k {k\choose i} p^i (1-p)^{k-i} x^i + o(1) \cdot \sum_{k=0}^\infty \sum_{i=0}^k {k\choose i} p^i (1-p)^{k-i} x^i\\
&=& p\sum_{k=0}^\infty a_k (1-p+px)^k + o(1) \sum_{k=0}^\infty (1-p+px)^k\\
&=& pA(1-p+px) + \frac{o(1)}{1 - x}.
\end{eqnarray*}\\

\noindent
{\bf Scale-free degree sequence starting at 0.}\\

A naive way to define scale-freeness is to require
\begin{equation}
  \label{eq:3}
  b_i= f(p)a_i + o(1) \quad \mbox{ for all } i\geq 0,
\end{equation}
where $f(p)$ is a quantity depending only on $p$.

Equivalently, for any $x\in [-1,1]$ and $p\in (0,1)$, we have
\begin{equation}
  \label{eq:4}
 pA(1-p+px)=f(p)A(x).
\end{equation}
To solve equation (\ref{eq:4}), let $x=1$.
We get $pA(1)=f(p)A(1)$. Thus $f(p)=p$. We have
\begin{equation}
  \label{eq:5}
  A(1-p+px)=A(x).
\end{equation}
Let $x=0$. We have $A(0)=A(1-p)$. Therefore,
\begin{eqnarray*}
  A'(0) &=& \lim_{x\to 0} \frac{A(x)-A(0)}{x}\\
&=& \lim_{x\to 0} \frac{A(1-p+px)-A(1-p)}{x}\\
&=&pA'(1-p).
\end{eqnarray*}
Since this holds for any $p\in (0,1)$,
we have
\begin{eqnarray*}
A(p) &=& A(0)+ \int_{1-p}^1 A'(1-p) \, dp  \\
&=& A(0) +  \int_{1-p}^1 \frac{A'(0)}{p} \, dp\\
&=& A(0) - A'(0)\ln(1-p).
\end{eqnarray*}
Thus, $$A(x)=A(0)-A'(0)\ln(1-x).$$
We have
\begin{eqnarray*}
A(1-p+px) &=& A(0) - A'(0)\ln(p-px)\\
&=& A(0) - A'(0)(\ln p + \ln (1-x))\\
&=& A(x) - A'(0)\ln p.
\end{eqnarray*}
This forces $A'(0)=0$. The only solution for equation (\ref{eq:4})
is $A(x) \equiv A(0)$ (the constant function, corresponding to a graph with no edges).  This solution is not interesting.\\

\noindent {\bf Scale-free degree sequence starting at 1.}\\

In many cases, we do not care about the number of isolated vertices.
We only require that
\begin{equation}
  \label{eq:6}
  b_d= f(p)a_d + o(1) \quad \mbox{ for all } d\geq 1.
\end{equation}
where $f(p)$ is a quantity depending only on $p$.

Equivalently, for any $p\in (0,1)$ and $x\in [-1,1]$, we have
\begin{equation}
  \label{eq:7}
  f(p)(A(x) -A(0)) = p(A(1-p+px)-A(1-p)).
\end{equation}
Take the derivative with respect to $x$ on both sides.
We have,  for any $p\in (0,1)$ and $x\in (-1,1)$,
\begin{equation}
  \label{eq:8}
  f(p)A'(x)=p^2A'(1-p+px).
\end{equation}
Let $\alpha=\int_0^1 \frac{f(p)}{p^2}dp$ be a positive constant.
Divide both sides of equation (\ref{eq:8}) by $p^2$ and integrate it
with respect to $p$ from $0$ to $1$.
We have
\begin{eqnarray*}
  \alpha A'(x) &=& \int_0^1 A'(1-p +px) dp \\
&=& \frac{A(1)-A(x)}{1-x}\\
&=& \frac{1-A(x)}{1-x}.
\end{eqnarray*}
Rewriting this expression,
\begin{equation}
  \label{eq:9}
\frac{A'(x)}{1-A(x)}= \frac{1}{\alpha(1-x)}.
\end{equation}
Now, integrate with respect to $x$ from $0$ to $x$. We get
\begin{equation}
  \label{eq:10}
 \ln \frac{1-A(0)}{1-A(x)} = - \frac{1}{\alpha} \ln (1-x).
\end{equation}
Therefore, we have
\begin{equation}
  \label{eq:11}
  A(x)=1- (1-A(0))(1-x)^\frac{1}{\alpha}.
\end{equation}
It is easy to verify that equation (\ref{eq:11}) satisfies
equation (\ref{eq:8}) with $f(p)=p^{\frac{1}{\alpha}+1}$.

We do not care about $A(0)=a_0$, the number of isolated vertices.  Hence,
the solution is uniquely determined by a parameter $\alpha$ up to a
a constant factor.
For $d\geq 1$, we have
\begin{eqnarray*}
a_d &=& (1-a_0) {\frac{1}{\alpha} \choose d} (-1)^{d+1}  \\
&=& -(1-a_0) {d-\frac{1}{\alpha} -1 \choose d}\\
&=& O(d^{-(\frac{1}{\alpha}+1)}).
\end{eqnarray*}
In other words, the degree frequency sequence follows a power-law distribution with
exponent $\beta=1+1/\alpha$.  However, not all $a_d$ are positive.  Particularly, if $\beta > 2$, then there are negative terms $a_d$, $d \geq 1$.\\

\noindent
{\bf Scale-free degree sequence starting at $k$.}\\

Now we assume that the degree sequence distribution, considering only degrees at least
 $k$, is scale-free.  That is,
\begin{equation}
  \label{eq:12}
  b_d= f(p)a_d + o(1) \quad \mbox{ for all } d\geq k.
\end{equation}
where $f(p)$ is a quantity depending only on $p$.

Or equivalently, for any $p\in (0,1)$ and $x\in [-1,1]$, we have
\begin{equation}
  \label{eq:13}
  f(p)(A(x) -\sum_{d=0}^{k-1} a_d x^d) = p(A(1-p+px)-\sum_{d=0}^{k-1} a_d
x^d{k\choose d}p^d(1-p)^{k-d}).
\end{equation}
Take the $k$-th derivative with respect to $x$ on both sides
to get rid of all terms of degree up to $k-1$.
We have,  for any $p\in (0,1)$ and $x\in (-1,1)$,
\begin{equation}
  \label{eq:14}
  f(p)A^{(k)}(x)=p^{k+1}A^{(k)}(1-p+px).
\end{equation}
Let $\alpha_k =\int_0^1 \frac{f(p)}{p^{k+1}} dp$.
Similar arguments to those above show that the solution of equation (\ref{eq:14}) is
of form
$$A^{k-1}(x)= C_1 - C_2 (1-x)^{\frac{1}{\alpha}}.$$
If we then integrate with respect to $x$ $k-1$ times, the result is
\begin{equation}
  \label{eq:15}
  A(x)=P_k(x) - C (1-x)^{\frac{1}{\alpha_k} + k}.
\end{equation}
Here $P_k(x)$ is a polynomial of $x$ with degree $k-1$.
It is easy to verify that equation (\ref{eq:15}) is the solution
of equation (\ref{eq:13}) with $f(p)=p^{\alpha_k+k}$.
Let $\beta=\frac{1}{\alpha_k}+k$.
For any $d\geq k$, we have
\begin{equation}
  \label{eq:16}
  a_d=C{d-\beta \choose d}
\end{equation}
If we set
$$
C = C_\beta = \left ( \sum_{d \geq \ceil{\beta}} \binom{d-\beta}{d} \right )^{-1}
$$
then the $a_d$ are positive for $d > \beta$.  Note that $\sgn(C_\beta) = (-1)^{\floor{\beta}}$.\\

\section{Concentration}

Since we know that the only degree sequences which are scale-free {\it in expectation} have power-law limit distributions, it is crucial to show that such graph have degree sequences which are close to their means with high probability.

\begin{theorem} Suppose that $\{G^n\}_{n=1}^\infty$ is a sequence of graphs on $n \rightarrow \infty$ vertices with degree sequence of limit distribution $\{\lambda_d\}_{d=k}^\infty$.  Further suppose that
$$
\sum_{v \in G} \deg(v)^2 = O(n^{2 - \epsilon})
$$
for some $\epsilon > 0$.  Then the degree sequence of $G^n_p$ also has a limit distribution $\{\lambda^\prime_d\}_{d=k}^\infty$.
\end{theorem}
\begin{proof}
Let $a_d = a_d(n)$ be the fraction of vertices of degree $d$ in $G^n$ and let $b_d = b_d(n)$ be the fraction of vertices of degree $d$ in $G^n_p$.  Clearly it suffices to show that $b_d$ is concentrated about its expectation.

To that end, we apply the Azuma-Hoeffding inequality to the ``vertex exposure'' martingale.  In particular, consider the following process.  Fix $d \geq k$, order the vertices of $G^n$ as $v_1,\ldots,v_n$, and let $A_m$ denote the event that $v_m \in G^n_p$.  Let $X_0 = \bE[b_d n]$, and let $X_{m+1} = \bE[X_m | A_{m+1}]$.  That is, at stage $m$, we ``expose'' vertex $m$ and recalculate the expected number of vertices of degree $d$ based on the new information concerning whether or not $v_m \in G^n_p$.  It is easy to see that this is a martingale, and, furthermore, that $|X_{m+1} - X_m| \leq \deg(v_{m+1}) + 1$, where $\deg(\cdot)$ denotes degree in $G$.  Since $b_d n = X_n$, we may apply the Azuma-Hoeffding inequality to get
$$
\bP \left [|b_d - \lambda_d| \geq t/n \right ] \leq \exp \left ( \frac{-t^2}{2 \sum_{m=1}^n (\deg(v_m) + 1)^2} \right )
$$
for $t \geq 0$.  Since $\sum_{m=1}^n \deg(v_m)^2 = O(n^{2-\epsilon})$ and 
$$
\sum_{m=1}^n \deg(v_m) \leq \sqrt{n} \left (\sum_{m=1}^n \deg(v_m)^2 \right)^{1/2} = O(n^{3/2 - \epsilon/2})
$$
by Cauchy-Schwarz, we can set $t = n^{1-\epsilon/4}$, getting
\begin{equation} \label{eq:17}
\bP \left [|b_d - \lambda_d| \geq t/n \right ] \leq e^{-\Omega(n^{\epsilon/2})}.
\end{equation}
Let $t^\prime = t/n = n^{-\epsilon/4}$.  Then, since
$$
\sum_{n = 1}^\infty \bP \left [\bigwedge_{d = k}^n (|b_d - \lambda_d| \geq t^\prime) \right ] \leq \sum_{n = 1}^\infty n e^{-n^{\epsilon/2}} < \infty,
$$
the Borel-Cantelli Lemma implies that asymptotically almost surely, $|b_d - \lambda_d| \leq t^\prime = o(n)$ for all $d \geq k$.

\end{proof}

We have the following theorem.
\begin{theorem}
  For any integer $k>\beta>1$, the degree sequence starting at $k$
  defined by $a_d=C_\beta {d-\beta \choose d}n+o(n)$ is
    scale-free.  Moreover, if a graph $G$ on $n$ vertices such that
$$
\sum_{v \in G} \deg(v)^2 = O(n^{2-\epsilon})
$$
for some $\epsilon > 0$ has a scale-free
    degree sequence starting at $k$, then there is a $\beta\in (1,k)$ so that $a_d=C_\beta {d-\beta \choose d}n
      +o(n)$. As a consequence, sparse graphs with scale-free degree sequence
are power-law graphs.
\end{theorem}

\section{Scale-free set system}
Many power-law graphs like the Collaboration Graph and the Hollywood
Graph are actually better modeled by set systems (or hypergraphs)
rather than graphs. For example, in the Math Reviews database, each
published item has one or more authors.  The family of all papers
considered as collections of authors forms a set system.  The Collaboration Graph only captures part of the
information in this set system.  Here we quote from the Erd\H{o}s
number project \cite{enp}:\\

\begin{minipage}[c]{0.9\textwidth} \it
There are about 1.9 million authored items in the Math Reviews
database, by a total of about 401,000 different authors. $\ldots$
Approximately
62.4\% of these items are by a single author, 27.4\% by two authors, 
8.0\% by three authors, 1.7\% by four authors, 0.4\% by five authors, 
and 0.1\% by six or more authors.
\end{minipage}\\

In this example, the distribution of set-sizes follows a power-law
distribution. Is this just a coincidence? Is ``scale-free''
distribution of a set system always a power-law distribution?
\begin{figure}[htbp]
  \centering
  \psfig{figure=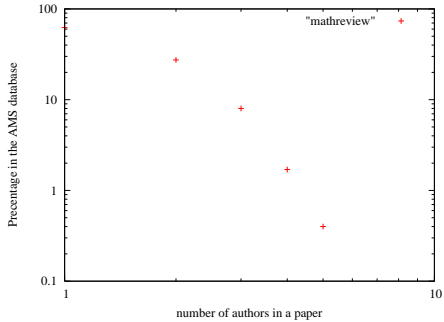, width=0.5\textwidth}
  \caption{The precentage of multiple-author-paper in AMS Review database.}
  \label{fig:mathreview}
\end{figure}

Motivated by this example and ``scale-free'' graphs,
we consider the following problem. 
For a set system ${\cal F}$ and any probability $p\in(0,1)$, 
the random sub-set-system ${\cal F}_p$ is 
chosen by independently removing vertices with probability $1-p$
 and reducing the sets to their remaining elements.

 \begin{problem}
   For what sequence of set-sizes in a set system ${\cal F}$,
is the sequence of the set-sizes in random sub-set-system ${\cal F}_p$
essentially the same as the original sequence up to a scale?
 \end{problem}

For $i\geq 1$, let $a_i$ be the number of $i$-sets in ${\cal F}$ 
and $b_i$ be the number of $i$-sets in ${\cal F}_p$. We are asking
if there is a function $f(p)$ such that
$$b_i=f(p)a_i+o(n)$$
for all $i\geq k$. Here $k$ is a small positive integer.

Since the expected value $\bE(b_i)$ satisfies
\begin{equation}
  \label{eq:18}
  \bE(b_i)=\sum_{j\geq i} a_j{j\choose i}p^i(1-p)^{j-i}.
\end{equation}

It is necessary to have
\begin{equation}
  \label{eq:19}
  \sum_{j\geq i} a_j{j\choose i}p^i(1-p)^{j-i}=f(p) a_i
\end{equation}
for all $i\geq k$.

Let $A(x)=\sum_{i}a_ix^i$ be the generating function.
For any $p\in (0,1)$ and $x\in [-1,1]$, we have
\begin{equation}
  \label{eq:20}
  f(p)(A(x) -\sum_{d=0}^{k-1} a_d x^d) = (A(1-p+px)-\sum_{d=0}^{k-1} a_d
x^d{k\choose d}p^d(1-p)^{k-d}).
\end{equation}
This is essentially the same equation as equation (\ref{eq:13}).
Thus we have the following theorem.

\begin{theorem}
If the sequence of set-sizes in a set-system starting at $k>1$
is scale-free, then there are constants $\beta\in (1,k)$ and $C$ such
that the number of $i$-sets in this set-system is $C_\beta {i-\beta \choose i}n+o(n)$ for all $i\geq k$.

\end{theorem}

\section{Remarks and questions}

Note that the results of the preceding sections have a probabilistic
interpretation.  Suppose that, for each $n$, we have a probability
distribution $\cG$ over graphs on $n$ vertices with the property that
the expected number of vertices of degree $d$ is $a_d$.  Then, what
must $\bE[a_d]$ be if, when $G$ is sampled from $\cG$ and a random
subgraph $G_p$ is taken, the expected number $b_d$ of vertices of
degree $d$ after scaling so that $\sum_d a_d = \sum_d b_d$ is the same
as $a_d$?  The above analysis provides the answer: the expectation of
$a_d$ must be a power law in $d$.

Now, it is natural to ask, if the variance of the $b_d$ is scaled as
the square of the scaling factor for the expectations, then what must
$\sigma^2(a_d)$ be?  In fact, one can ask the same question of all
moments, leading to the following open problem:

\begin{problem} Fix $p \in (0,1)$.  Let $G$ be drawn from a
  probability distribution $\cG$ on graphs with $n$ vertices.  Suppose
  that $a_d$, $d \geq 0$, is the number of vertices of degree $d$ in
  $G$, and $b_d$, $d \geq 0$, is the number of vertices of degree $d$
  in $G_p$.  For which distributions $\cG$ is it true that there
  exists some $c(p) \in \mathbb{R}$ so that $\{a_d\}_{d \geq k}$ and
  $\{c(p) b_d\}_{d \geq k}$ have approximately the same distribution
  for large $n$?  Is it possible to find such $\cG$ for all $p \in
  (0,1)$ simultaneously?
\end{problem}

Currently, the exponents of ``real-world'' scale-free networks' power
laws is estimated in a rather ad-hoc fashion, usually using a
regression on the log-log plot of frequency vs.\@ degree after
removing the extremes of the data.  If it were possible to describe
scale-free distributions exactly, then it would make sense to ask the
following very practical question:

\begin{problem} Find an unbiased estimator for the exponent of a power-law degree distribution.
\end{problem}

For the matter of the variance of the $a_d$, we note that, at least
for $\beta \in (1,2)$, the following must be true: $$
p^{2\beta}
\sigma^2(a_d) = \sum_{k} \binom{k}{d}^2 p^{2d} (1-p)^{2k-2d} \left
  (\sigma^2(a_k) + \binom{k-\beta}{k} \right) - p^{\beta-1}
\binom{d-\beta}{d}.  $$
This statement can be proven by applying the
formula $$
\sigma^2(\sum_{i=1}^N X_i) = \bE[X_1]^2 \sigma^2(N) +
\bE[N] \sigma^2(X_1) $$
for i.i.d.\@ variables $X_i$ and an
independent variable $N$ taking on nonnegative integer values.

We also ask, what can be proved by extending the definition of
scale-freeness to hypergraphs?  We believe that the situation is very
similar to that of graphs when the hypergraphs being considered are
uniform (with edges removed whenever at least one of their vertices is
removed).  Perhaps the answer lies in a more refined description of scale-freeness.  For example, consider the quantity $a_H(G)$, the number of occurrences of $H$ as an induced subgraph of $G$.  Suppose that $a_H(G)/n \rightarrow \alpha_H$ for each $H$ and some $\alpha_H \in \mathbb{R}^+$, and that this sequence is scale-free, i.e.,
$$
a_H(G_p) \propto a_H(G)
$$
for any fixed $p$ with $0 < p < 1$ and $H$ varying over all graphs on at least $k$ vertices.  Then what must $G$ look like?


\begin{thebibliography}{99}

\bibitem{acl}
W. Aiello, F. Chung and L. Lu,
A random graph model for massive graphs,
{\em Proceedings of the Thirty-Second Annual ACM Symposium on
  Theory of Computing}, (2000) 171-180.

\bibitem{aclp}
W. Aiello, F. Chung and L. Lu,
A random graph model for power law graphs,
{\it Experimental Math.}, {\bf 10}, (2001), 53-66.


\bibitem{acle}
W. Aiello, F. Chung and L. Lu,
Random evolution in massive graphs,
{\it Handbook of Massive Data Sets}, Volume 2, (Eds. J. Abello et al.),
 Kluwer Academic
Publishers, (2002), 97-122.

\bibitem{as} 
N. Alon, J. H. Spencer, 
The probabilistic method.
Second edition.
Wiley-Interscience [John Wiley \& Sons], New York,  2000.

\bibitem{ba}
A.-L.  Barab\'asi and R.  Albert,
Emergence of scaling in random networks,
{\rm Science} {\bf 286} (1999) 509-512.

\bibitem{baj}
A.~Barab\'{a}si, R.~Albert, and H.~Jeong,
Scale-free characteristics of random networks:
the topology of the world wide web,
{\it Physica} a {\bf 281} (2000), 69-77.

\bibitem{broder} A.~Broder, R.~Kumar, F.~Maghoul, P.~Raghavan,
S.~Rajagopalan, R.~Stata, A.~Tompkins, and J.~Wiener,
``Graph Structure in the Web,'' {\em proceedings of the WWW9 Conference},
May, 2000, Amsterdam.

\bibitem{br} B. Bollab\'as and O. Riordan,
Robustness and Vulnerability of Scale-free Random Graphs.
{\em Internet Mathematics}, {\bf 1} (2003) no. 1, 1--35.
\bibitem{ch0}
F. Chung, {\it Spectral Graph Theory},
 AMS Publications, 1997.

\bibitem{cl1}
F. Chung and L. Lu,
The diameter of random sparse graphs,
{\it Advances in Applied Math.} {\bf 26} (2001), 257-279.

\bibitem{conn}
F. Chung and L. Lu,
Connected components in a random graph with given degree
sequences, {\it Annals of Combinatorics} {\bf 6}, (2002), 125--145.

\bibitem{coupling}
F. Chung and L. Lu,
Coupling on-line and off-line analyses for
random power law graphs,
{\it Internet Mathematics}, to appear.

\bibitem{ave}
F. Chung and L. Lu,
The average distance in random graphs with given expected degrees,
{\it Proceedings of the National Academy of Science} {\bf 99}, (2002),
15879--15882.

\bibitem{spectra}
F. Chung, L. Lu and V. Vu,
The spectra of random graphs with given expected degrees,
{\it Proceedings of National Academy of Sciences}, {\bf 100}, no. 11, (2003), 63
13-6318.

\bibitem{bio} F. Chung, L. Lu, T. G. Dewey, and D. J. Galas,
Duplication models for biological networks,
{\it Journal of Computational Biology}, {\bf 10}, No. 5, (2003), 677-688.


\bibitem{cf}
C. Cooper and A. Frieze,
A general model of undirected Web graphs, {\it Random Structures and Algorithms}, {\bf 22},
(2003), 311-335.

\bibitem{cfv}
C. Cooper, A. Frieze and J. Vera,
Random vertex deletion in a scale free random graph, preprint.

\bibitem{eg} P. Erd\H{o}s and  T. Gallai,
Gr\'afok el\H{o}\'irt fok\'u pontokkal (Graphs with points
of prescribed degrees, in Hungarian), {\it Mat. Lapok} {\bf 11} (1961),
264-274.\bibitem{er} P. Erd\H{o}s and A. R\'enyi, On random graphs. I,
{\em
    Publ. Math. Debrecen} {\bf 6} (1959), 290-291.
    
\bibitem{enp} J.\@ Grossman, P.\@ Ion, and R.\@ De Castro, The Erd\H{o}s Number Project,
\texttt{http://www.oakland.edu/enp/}.

\bibitem{jlr}
S. Janson, T. {\L}uczak, and A. Rucinski,
{\it Random Graphs}, Wiley-Interscience, 2000.

\bibitem{jta}
H. Jeong, B. Tomber, R. Albert, Z. Oltvai and A. L. Bab\'arasi,
The large-scale organization of metabolic networks,
{\it Nature}, {\bf 407} (2000), 378-382.

\bibitem{klein}
J. Kleinberg, S. R. Kumar, P. Raphavan, S. Rajagopalan and A. Tomkins,
The web as a graph: Measurements, models and methods,
{\it Proceedings of the International Conference on Combinatorics
and Computing}, 1999.

\bibitem{kv}
J.H. Kim
and
V. Vu,
 Concentration of multi-variate
  polynomials and its applications, {\em Combinatorica}, {\bf 20} (3) (2000),
  417-434.

\bibitem{lotka}
A. J. Lotka,
The frequency distribution of scientific productivity,
{\it The Journal of the Washington Academy of the Sciences},
{\bf 16} (1926), 317.

\bibitem{lu}
 L. Lu, The Diameter of Random Massive Graphs,
{\em Proceedings of the Twelfth ACM-SIAM Symposium on
Discrete Algorithms}, (2001) 912-921.

\bibitem{mcdiarmid}
 C. McDiarmid,
 Concentration,
{\em Probabilistic methods for algorithmic discrete mathematics},
 195--248, {\em Algorithms Combin.}, 16, Springer, Berlin, 1998.

\bibitem{milgram}
S. Milgram,
The small world problem,
{\it Psychology Today}, {\bf 2} (1967), 60-67.
\bibitem{mit}
M. Mitzenmacher,
A brief history
of generative models for power law and lognormal distributions,
{\it Internet Mathematics}, to appear.

\bibitem{stumpf}
M.\@ P.\@ H.\@ Stumpf, C.\@ Wiuf, R.\@ M.\@ May,
Subnets of scale-free networks are not scale-free: Sampling properties of networks,
{\it Proceedings of the National Academy of Sciences},
{\bf 102} (2005), 4221--4224.

\end{thebibliography}
\end{document}